\documentclass[reqno]{amsart}
\usepackage{amssymb}
\usepackage{graphicx}
\usepackage{xcolor} 
\usepackage{fullpage} 
\usepackage{amsmath}
\usepackage{amsthm}
\usepackage{verbatim}
\usepackage{enumitem}
\usepackage{mathrsfs}

\usepackage{mathtools}
\mathtoolsset{showonlyrefs} 

\graphicspath{{Fig/}}

\newtheorem{thm}{Theorem}

\newtheorem{prop}[thm]{Proposition}

\newtheorem{fact}[thm]{Fact}

\theoremstyle{definition}
\newtheorem{dfn}{Definition}
\newtheorem{rem}{Remark}
\newtheorem{example}{Example}

\DeclareMathOperator{\dvg}{div}
\DeclareMathOperator{\grad}{grad}
\DeclareMathOperator{\jgrad}{\mathcal{J}grad}

\newcommand{\vf}{\mathfrak{X}}
\newcommand{\Diff}{\mathrm{Diff}}

\newcommand{\Fix}{\mathrm{Fix}}

\newcommand{\rp}{\mathrm{Re}}

\newcommand{\re}{{\mathbb R}}
\newcommand{\ce}{{\mathbb C}}
\newcommand{\ze}{{\mathbb Z}}
\newcommand{\te}{{\mathbb T}}
\newcommand{\se}{{\mathbb S}}
\newcommand{\me}{\mathrm{e}}
\newcommand{\mi}{\mathrm{i}}
\newcommand{\md}{\mathrm{d}}

\newcommand{\tc}{\Gamma^*}

\newcommand{\lie}{\mathcal{L}}
\newcommand{\rie}{\mathcal{R}}

\begin{document}

\title{Point vortex on surfaces with continuous symmetry}

\author{Yuuki Shimizu} 
\address{Faculty of Science, Academic Assembly,
University of Toyama, 3190 Gofuku, Toyama 930-8555, Japan} 
\email{shimizu@sci.u-toyama.ac.jp}


\date{\today} 


\begin{abstract}
We derive an analytic formula for the hydrodynamic Green function and the Robin function on every orientable surface admitting a hydrodynamic Killing vector field. Closed-form expressions are provided for all fourteen canonical Riemann surfaces, covering both compact and non-compact cases; the formulae satisfy the slip boundary condition and generate complete Hamiltonian vector fields.

As an application, we clarify the mechanism whereby the curvature affects a point vortex in both qualitative and quantitative viewpoints. 
Qualitatively, we show a single point vortex is governed by a Hamiltonian flow whose vorticity is given by the curvature up to area constant. 
Quantitatively, on a rectangular torus with periodic curvature we use the analytic formula to describe two regimes: linear response that mirrors the curvature wave when the mean component is small, and a nonlinear response with amplitude resonance. 
The results supply a unified tool for detailed studies of point vortex dynamics and Euler–Arnold flows on surfaces with continuous symmetry.
\end{abstract}

\maketitle

\section{Introduction}
Explicit Green functions for the Laplace operator on curved two-dimensional manifolds are indispensable well beyond pure analysis: they enter boundary-integral methods in incompressible fluid mechanics, potential and electromagnetic theory on curved substrates, and numerical solvers for thin-film models in materials science. Whenever the underlying surface deviates from the plane, having a closed analytic kernel immediately improves accuracy, reduces computational cost, and often reveals hidden geometric structure. Yet on most non-trivial surfaces only fragmentary formulas are available, and many are hard-coded for a single application or boundary condition.

In two–dimensional incompressible and inviscid fluid dynamics on a curved surface, a particular variant—the \emph{hydrodynamic Green function}(HGF)—is central, because every velocity field may be written as the convolution of the vorticity with this kernel. Detailed information on the Green function therefore translates directly into quantitative control of the flow. 
A classical way to exploit that control is to approximate a smooth vorticity distribution by a finite sum of Dirac measures, the \emph{point vortices}.  Advected by the Euler–Arnold equation, these vortices move according to a finite–dimensional Hamiltonian system whose Hamiltonian is expressed through the HGF and its regularization by geodesic distance.  Level sets of the Green function coincide with particle paths in the presence of a single point vortex, so explicit or computable formulae have long been sought.  Such formulae are known on selected surfaces: spheres and hyperbolic discs~\cite{Kimura_Okamoto_1987, Kimura_1999}, multiply–connected plane domains~\cite{Sakajo_2009}, cylinders~\cite{Montaldi_Souliere_Tokieda_2003}, flat tori~\cite{Tkachenko_1966}, the Bolza surface~\cite{Ragazzo_2017}, surfaces of revolution~\cite{Hally_1980,Dritschel_Boatto_2016}, and toroidal surfaces~\cite{Sakajo_Shimizu_2016}, among others.  They underpin models of geophysical flows~\cite{Nycander_1996,Nasser_Sakajo_Murid_Wei_2015} and quantized vortices in superfluid films~\cite{Corrada-Emmanuel_1994}.

Despite these examples, a systematic understanding of how surface geometry constrains vortex dynamics is still incomplete.  Some phenomena appear universal: two–vortex motion shows similar qualitative behavior on the plane, the sphere, and hyperbolic space~\cite{Kimura_1999}, and polygonal $N$–ring configurations share stability properties across several constant–curvature geometries~\cite{Montaldi_Tokieda_2013}. Yet even closely related families can differ: on a one–parameter family of tori the linear stability of an $N$–ring changes qualitatively~\cite{Sakajo_Shimizu_2018}.  A unified perspective is needed to decide which effects stem from global geometry and which are surface–independent.

One pragmatic hypothesis in recent work is to assume the flow field admits a non–trivial Killing vector field.  Such vector fields preserve the Riemannian metric and are stationary solutions of the Euler–Arnold equation.  They generalize uniform flows and zonal jets, and their existence ensures slip boundary conditions and completeness of trajectories. 
We focus on surfaces that admit a hydrodynamic Killing vector field (HKVF), which is a nonzero complete Killing vector field satisfying the slip boundary condition.
Every orientable surface admitting an HKVF is conformally equivalent to one of the 14 canonical Riemann surfaces equipped with a rotationally or translationally symmetric metric \cite{Shimizu_2023_Killing}.
Thus symmetry—rotational or translational—is the sole geometric constraint: once the symmetry is fixed, the metric may vary freely, yet the presence of an HKVF guarantees an underlying continuous isometry that is invaluable for analyzing fluid motion.

The present paper continues that program.  Our first main result, Theorem~\ref{thm-Analytic-form}, provides analytic formulae of the hydrodynamic Green function $G$ and the Robin function $R$ on every surface admitting an HKVF.  
The formula unifies all previously known closed cases and extends them to non-compact surfaces with ends, while automatically fulfilling the slip condition and completeness requirements.

As an application of Theorem \ref{thm-Analytic-form}, Section \ref{sec-impact-curvature} clarifies the mechanism whereby the curvature affects a point vortex in both qualitative and quantitative viewpoints. The key identity proved in \ref{subsec-general-curvature},
\begin{align}
	-\triangle R=\frac{\kappa}{2\pi}+2c_{M},
\end{align}
where $\kappa$ is the Gauss curvature and $c_{M}$ is the area constant. Hence, the motion of a single point vortex is identified with a Hamiltonian flow whose vorticity is given by the curvature up to area constant. 

Building on the qualitative observations from the previous section, we now quantitatively examine the effect of curvature on the dynamics of a point vortex in greater detail.  
Choosing the torus \(M=\mathbb{T}_{\Lambda}\) with rectangular lattice \(\Lambda=\{\,2\pi\mi m+Tn\,\}\) and adopting the periodic profile \(\kappa(x,y)=\kappa_{0}+\kappa_{1}\cos\omega x\), the analytic formulas supplied by Theorem~\ref{thm-Analytic-form} yield a closed expression for the vortex velocity.  
When \(\kappa_{0}\to0\) the asymptotic formula
\[
v=\frac{\kappa_{1}}{2\pi\omega^{2}}\!
\left(\frac{\kappa_{1}}{\kappa_{0}}+2\cos\omega x+\Bigl(\omega-\frac{4}{T}\Bigr)\sin\omega x\right)
\]
reveals a \emph{linear response}: the output velocity shares the same periodicity as the input curvature, and an \emph{affine similarity law} links their waveforms with an explicit similarity ratio and phase shift.  
Near the resonance \(|\kappa_{1}|-|\kappa_{0}|\ll1\) we instead find
\[
v=-\frac{K}{2\pi\omega^{2}}
\left(-\omega P_{1}(x)+A(x)P_{2}(x)\right),
\]
where the double--angle terms in \(P_{1},P_{2}\) generate a distinct \emph{waveform asymmetry}.  
The competition between \(\omega P_{1}\) and \(A P_{2}\) explains how the period \(T\) and lattice parameter \(N\) switch the dominant mode, completing a quantitative picture of the impact of the curvature on a point vortex.


The paper is organized as follows.  Section~\ref{sec-prel} recalls basic facts about surfaces, Killing fields, and Euler–Arnold flows. Section~\ref{sec-main} states and proves the analytic formula for $G$ and $R$.  Section~\ref{sec-impact-curvature} establishes the curvature–vortex relation by analyzing the linear and nonlinear regimes on the torus.  Section~\ref{sec-conclusion} summarizes the results and outlines further directions.

\section{Preliminaries}
\label{sec-prel}
We call a connected orientable Riemannian $2$-manifold $(M,g)$ possibly with boundary as a surface in this paper. 
Conventionally, for boundary-less surfaces, the surfaces are called closed surfaces when compact, and open surfaces otherwise. 
For non-compact surfaces, unbounded connected components of the complement of a compact subset in the surface are called ends of the surface

Let $\mathfrak{X}(M)$ be the space of all vector fields $X:M\to TM$. 
Let $\mathrm{Diff}(M)$ be the space of all diffeomorphisms on $M$. 
The time $t$-map of each vector field $X$ is denoted by $X_t\in\mathrm{Diff}(M)$. 
A vector field $X\in\mathfrak{X}(M)$ is said to be complete if for each $t\in\re$, $X_t\in\mathrm{Diff}(M)$ is defined. 
Let $\mathrm{Fix}(X)=\{p\in M|\, X(p)=0\}$ and $\mathrm{Per}(X)=\{p\in M|\,\exists  T \in (0,\infty),\, X_T(p)=p\}$. 
The slip boundary condition for a vector field $X$ on the boundary $\partial M$ is that $g(X,n)=0$ on $\partial M$ where $n$ is the inward unit normal vector for $\partial M$, or equivalently, $X|_{\partial M}\in\mathfrak{X}(\partial M)$. 
An orientation preserving diffeomorphism $\varphi: M\to N$ between surfaces $(M,g)$ and $(N,h)$ is called a conformal mapping if there exists a positive function $\lambda$ on $M$ such that 
\begin{align}
	\varphi^*h =\lambda^2 g. 
\end{align}
Two surfaces $(M, g)$ and $(N, h)$ are said to be conformally equivalent, denoted as $M \simeq N$, if there exists a conformal mapping $\varphi: M \to N$. 

Note that 
\begin{align}
\label{eq-normal-G}
	*\md G=-\mi \partial G\md z+\mi \bar\partial G\md \bar z=2\rp(-\mi \partial G\md z).  
\end{align}
\begin{align}
	\jgrad G=\lambda^{-2}(-2\mi \partial G\bar \partial+2\mi \bar \partial G\partial) =4\lambda^{-2}\rp( -\mi \partial G\bar\partial). 
\end{align}

%


\section{Euler-Arnold flows}
\subsection{Euler-Arnold equation on surfaces}
\label{sec-EA}
The Euler-Arnold equation describes the motion of incompressible and inviscid fluid flows: 
\begin{align}
\label{eq-EA}
	v_t+\nabla_v v=-\grad p,\quad \dvg v=0         
\end{align}
subject to the slip boundary condition if $M$ has boundary. 
The velocity $v: [0,T]\times M\to TM$ and the pressure $p:[0,T]\times M\to \re$ are unknown except for their initial datum. 
Let us rewrite Eq.~\eqref{eq-EA} with respect to the flat operator $v^\flat= g(v,\cdotp)$. 
The Levi-Civita connection $\nabla$ has the following decompositions.
\begin{align}
\label{eq-decomp-nabla}
	\nabla_v v^\flat 
	&= \lie_vv^\flat -\md |v^\flat |^2/2.
\end{align}
Using the decomposition~\eqref{eq-decomp-nabla}, the Euler-Arnold equation~\eqref{eq-EA} is written as 
\begin{align}
	v_t^\flat+ \lie_v v^\flat -\md |v^\flat |^2/2
	&=-\md p,\label{eq-EA-sym}.
\end{align}

By Eq.~\eqref{eq-EA-sym} and $\lie_Xg =0$, 
for each Killing vector field $X\in\vf(M)$, $(X,p_X)$ with $p_X=|X|^2/2$ is a steady solution for the Euler-Arnold equation.
Owing to $\dim M=2$, $\omega= *\md v^\flat$ becomes a scaler function, called the vorticity.
Applying $*\md$ to Eq.~\eqref{eq-EA-sym}, we obtain the vorticity equation,
\begin{align}
\label{eq-vort}
	\omega_t+\lie_v\omega=0.
\end{align}
Hence, the Euler-Arnold equation on surfaces is regarded as the advection equation for the vorticity. 
To determine the characteristic curve, we conventionally derive the velocity $v$ from the vorticity $\omega$ by the Biot-Savart law. 
The Biot-Savart law is based on the fact that every incompressible vector field $v$ is a Hamiltonian vector field, that is, there exists a function, called a stream function, such that 
\begin{align}
	v&=-\jgrad \psi,\\
	-\triangle \psi &=\omega
\end{align}
as long as $M$ is a simply connected surface. 

Otherwise, the above approach does not provide general solutions for the Euler-Arnold equation.  
The problem arises even if the flow field is not curved, for instance, multiply connected domains in the plane.
Recently, some treatises for non-Hamiltonian solutions have been developed~\cite{Iftimie_Filho_Lopes2020,Gustafsson_2022,Ragazzo_2022,Shimizu_2023}. 
Of course, it is significant to focus on Hamiltonian solutions when exploring hydrodynamics from the viewpoint of special solutions. 
Therefore, in the present paper, we only consider a Hamiltonian solution as a special solution for the Euler-Arnold equation.

\section{Single point vortex}
If $v$ is given by $v=-\jgrad \psi$ for some stream function $\psi$, then 
\begin{align}
	\omega = *\md v^\flat=*\md (-*\md \psi)= -\triangle \psi. 
\end{align}
Using a Green function $G$ for the Laplacian, $\psi= \langle G,\omega\rangle$ solves this Poisson equation. 
Hence, $v$ is explicitly written by $\omega$ with certain regularity as 
\begin{align}
	v=-\jgrad \langle G,\omega\rangle,
\end{align}
which is a generalization of the Biot-Savart law to the curved surface. 
In particular, formally injecting a delta function $\delta_{x_0}$ centered at $x_0\in M$ into $\omega$, we derive a formal velocity field 
\begin{align}
	v(x)= -\jgrad_x G(x,x_0).
\end{align}
Then, the vorticity $\omega=\delta_{x_0}$ is called a point vortex with unit circulation. 
Since the Green function $G$ has the $\log$ singularity at $x_0$, i.e., $G(x,x_0)\approx (-1/2\pi)\log d(x,x_0)$, the vector field $v$ does not defined at $x_0$. 
By introducing the Robin function $R$
\begin{align}
	R(x)=\lim_{x_0\to x}\left(G(x,x_0)+\frac{1}{2\pi} \log d(x,x_0)\right) 
\end{align}
as a regularization of the Green function, we consider that the motion of the single point vortex is governed by 
\begin{align}
\label{eq-jgradR}
	\dot x=-\jgrad R(x).
\end{align}


\subsection{Hydrodynamic Green function}
We have seen that for Hamiltonian solutions, a Green function for the Laplacian is used to recover the velocity  from the vorticity. 
Especially, it is reasonable that the velocity $-\jgrad_x G(x,x_0)$ for a point vortex satisfies the slip boundary condition and becomes a complete vector field for each time to be a solution for the Euler-Arnold equation in certain sense. 
From this reason, a special Green function, called a hydrodynamic Green function, is well-investigated in hydrodynamics on surfaces.
\begin{dfn}[Hydrodynamic Green function]
\label{def-HGF}
	Let $(M,g)$ be a surface. 
	A function $G\in C^\infty(M\times M\setminus \{x=x_0\})$ is called a \textit{hydrodynamic Green function}(shortly, HGF) if the following conditions are satisfied: for each $(x,x_0)\in M\times M\setminus \{x=x_0\}$, 
	\begin{enumerate}
		\item $-\triangle G(x,x_0)=\delta(x,x_0)+c_M$.
		\item $G(x,x_0)=G(x_0,x)$.
		\item $-\jgrad G(\cdotp,x_0)|_{\partial M}\in\vf(\partial M)$. 
		\item For any $t\in\re$, $(-\jgrad G(\cdotp,x_0))_t\in\Diff(M\setminus\{x_0\})$ and if $x_0\not\in\partial M$, $(-\jgrad G(\cdotp,x_0)|_{\partial M})_t\in\Diff(\partial M)$
	\end{enumerate}
	where 
		\begin{align}
		c_M
		&= 
		\begin{dcases}
		-|M|^{-1},\quad &\text{if $M$ is a closed surface with the area $|M|$},\\
		0,\quad &\text{otherwise.}
		\end{dcases}
		\end{align} 
\end{dfn}


\begin{rem}
\label{rem-HGF}
We shall see that the definition stated in the above is a generalization of the hydrodynamic Green function in the preceding works~\cite{Aubin_1998,Flucher_1999,Ragazzo_Viglioni_2017}. 
Obviously, for closed surfaces, it is exactly the Green function for the Laplacian~\cite{Aubin_1998}. 
Concerning compact surfaces with boundary, let us recall that each of orbits of $-\jgrad G(\cdotp,x_0)$ is a connected component of contour lines $\{x\in M\setminus \{x_0\}|\, G(x,x_0)=c\}$ since it is a Hamiltonian system. 
Then, the third condition in Definition~\ref{def-HGF} is also written as the Dirichlet constant value condition, 
\begin{align}
	G(x,x_0)\equiv c(x_0,\gamma)
\end{align}
for some constant $c$ depending on $x_0$ and each connected component $\gamma$ of the boundary. 
Then, the value of these constants are determined to satisfy the Neumann condition, which is derived from $-\triangle G=\delta_{x_0}$, that is, 
\begin{align}
	1
	=\int_M -*\triangle G
	=\int_M -\md *\md G
	=\int_{\partial M} -*\md G. 
\end{align}

For each connected component $\gamma$ of $\partial M$, the quantity $\int_\gamma*\md G$ is called the circulation around $\gamma$, which satisfies  
\begin{align}
	\int_\gamma*\md G= \int_\gamma \partial_n G \md s. 
\end{align}
Hence, the total circulation $\tc$ given as 
\begin{align}
	\tc =\sum_{\gamma\subset\partial M}\int_\gamma *\md G
\end{align}
satisfies 
\begin{align}
	\tc =-1, 
\end{align}
which corresponds to the definition for the hydrodynamic Green function given in~\cite{Flucher_1999}. 
Since every smooth vector field on a compact manifold is complete and a finite level set $L(K)=\{x\in M|\, |G(x,x_0)| \leq K\}\subset M$ is a compact submanifold for sufficiently large $K>0$, $-\jgrad G(\cdotp,x_0)|_{L(K)}$ is a complete vector field on $L(K)$. 
Especially, applying this to the restricted vector field on $\partial M$, which is compact in $M$, $-\jgrad G(\cdotp,x_0)|_{\partial M}$ is complete on $\partial M$. 
Obviously, $|G(x,x_0)|\to \infty$ as $x\to x_0$, which yields we can take a compact exhaustion $L(K_1)\subset L(K_2)\subset \ldots $ with a monotone increasing sequence $K_n\nearrow\infty$ such that $M\setminus\{x_0\} =\bigcup_n L(K_n)$.
From this, we obtain the completeness of $-\jgrad G(\cdotp, x_0)$ on $M\setminus\{x_0\}$. 

For non-compact surfaces, the notion of HGF has recently established by Ragazzo and Viglioni~\cite{Ragazzo_Viglioni_2017} under the assumption that the target surfaces are of topological finite type (i.e. finite genus and finitely many ends) and only with parabolic or hyperbolic ends. 
Let us review their definition and confirm that our definition becomes its generalization. 
An end is said to be parabolic (resp. hyperbolic) if the end is conformally equivalent to the punctured disc $\Delta^*$ (resp. an open annulus $\Delta_\rho$). 

Then, there exists a function $G$ smoothly defined on $M\times M\setminus\Delta$ such that $G$ satisfies the first and second condition in Definition~\ref{def-HGF} and for each parabolic end $\gamma$, 
\begin{align}
\label{eq-para-end}
	\lim_{|z|\to 0} \left( G(z,x_0)-\frac{\Gamma_p}{2\pi}\log |z|\right)=c
\end{align}
and for each hyperbolic end $\gamma$, 
\begin{align}
\label{eq-hyp-end}
	\lim_{|z|\to \rho}G(z,x_0)=c,\quad \lim_{r\to\rho} \int_{|z|=r} *\md G =\Gamma_h
\end{align}
where $c=c(x_0,\gamma)$ is a constant and $\Gamma$ is the prescribed circulation, 
subject the total circulation equals $-1$, including the ends as well as the boundaries. 

What must be demonstrated here is that $-\jgrad G(x,x_0)$ is a complete vector field when $G$ satisfies Eq.~\eqref{eq-para-end} or Eq.~\eqref{eq-hyp-end}. 
If $M$ has a parabolic end $\gamma$, setting 
\begin{align}
	h(z)=G(z,x_0)-\frac{\Gamma}{2\pi }\log|z|,
\end{align}
we deduce from Eq.~\eqref{eq-para-end} that $h$ is bounded. 
Since $h$ is harmonic except for $z=0$, $h$ has a removable singularity at $z=0$.
Hence, taking the smooth extension $\tilde h$ of $h$ to $z=0$, we see that $G(z,x_0)=\tilde h(z)+(\Gamma/2\pi) \log|z| $ is smoothly defined in the chart except for $z=0$ and $G(z,x_0)\to \infty$ as $|z|\to \infty$. 
In the same way as the compact case, we deduce the completeness of $-\jgrad G(\cdotp, x_0)$ by taking the same compact exhaustion $L(K_n)$. 
When $M$ has a hyperbolic end $\gamma$, Eq.~\eqref{eq-hyp-end} gives a smooth extension $\tilde G$ of $G$ to $\{|z|=\rho\}$. 
It follows from $\tilde G|_{\{|z|=\rho\}}=c$ that $-\jgrad \tilde G(\cdotp,x_0) $ is complete, which gives the completeness of $-\jgrad G(\cdotp,x_0)$.

Consequently, for non-compact surfaces with only parabolic ends or hyperbolic ends, it can be deduced that a function $G$ satisfying the conditions (1) and (2) in Definition~\ref{def-HGF} also meets the conditions (3) and (4) if it satisfies Eq.~\ref{eq-para-end} for parabolic ends or Eq.~\ref{eq-hyp-end} for hyperbolic ends, respectively. 
Since almost all the surface targeted in this paper are included in this category, we will show that the function $G$ satisfies Eq.~\eqref{eq-para-end} or Eq.~\eqref{eq-hyp-end}, thereby meeting the third and fourth conditions in Definition~\ref{def-HGF}. 
Nevertheless, as the channel domains, for instance, $[0,2\pi)\times\re$ or $[0,2\pi]\times \re$ are not included in that category, we will directly show that the function $G$ fulfills the third and fourth conditions in Definition~\ref{def-HGF}, which serves as the rationale behind proposing a slightly more generalized definition.
On what surfaces HGF exists, see~\cite{Aubin_1998,Sario_Nakai_1970,Ragazzo_Viglioni_2017} and the references there in. 
\end{rem}


\section{Hydrodynamic Killing vector fields}
\label{sec-hkvf}
In the present paper, as flow fields, we focus on surfaces admitting a hydrodynamic Killing vector field. 
The notion of hydrodynamic Killing vector field is established in~\cite{Shimizu_2023_Killing} and introduced as an appropriate Killing vector field in the context of fluid dynamics. 
The Killing vector field $X$ is defined as a vector field preserving the Riemannian metric $g$, that is, $\mathcal{L}_X g=0$. 
As we see in Section~\ref{sec-EA}, Killing vector fields become steady solution for the Euler-Arnold equation if  they satisfy the slip boundary condition. 
In addition, in order to avoid any fluid particle leaking out of the flow field in finite time, the Killing vector fields should be complete vector fields. 
\begin{dfn}[Hydrodynamic Killing vector field]
	Let $(M,g)$ be a surface. 
	$X\in\mathfrak{X}(M)$ is called \textit{a hydrodynamic Killing vector field} (HKVF for short) if the following conditions are satisfied.
	\begin{enumerate}
		\item $\mathcal{L}_X g=0$.
		\item $X\not\equiv 0$. 
		\item $X|_{\partial M} \in\mathfrak{X}(\partial M)$.
		\item For any $t\in\re$, $X_t\in\mathrm{Diff}(M)$ and $(X|_{\partial M})_t \in\mathrm{Diff}(\partial M)$. 
	\end{enumerate}
\end{dfn}

\begin{example}
\label{ex-surf}
	Let $\se^1=\re/2\pi\ze$. 
	The following Riemann surfaces with a flat metric $g$, which gives constant curvature admit a hydrodynamic Killing vector field. 
	\begin{enumerate}
		\item The Riemann sphere $\hat\ce=\ce\cup\{\infty\}$.
		\item The complex plane $\ce$.
		\item The unit open disc $\Delta=\{z\in\ce |\, |z|<1\}$, or equivalently, the right half plane $(0,\infty)\times \re$ or the open channel $(0,2\pi)\times \re$.
		\item The punctured plane $\ce^*=\ce\setminus\{0\}$, or equivalently, the cylinder $\se^1\times\re$.
		\item The punctured open disc $\Delta^*=\Delta\setminus\{0\}$. 
		\item An open annulus $\Delta_\rho=\{z\in\ce |\, \rho<|z|<1\}$, $\rho \in(0,1)$. 
		\item A torus $\te_\Lambda=\ce/\Lambda $, $\Lambda =\{m\pi_1+n\pi_2|\,\pi_1,\pi_2\in\ce,\,\mathrm{Im}(\pi_1/\pi_2)>0,\,m,n\in\ze \}$.
		\item The unit closed disc $\overline\Delta=\{z\in\ce|\, |z|\le1\}$ 
		\item The punctured closed disc $\overline\Delta^*=\overline\Delta\setminus\{0\}$. 
		\item A closed annulus $\overline\Delta_\rho=\{z\in\ce|\, \rho\le |z|\le 1\}$. 
		\item A semi-closed annulus $\tilde \Delta_\rho=\{z\in\ce|\, \rho< |z|\le 1\}$.
		\item The closed right half plane $[0,\infty)\times\re$. 
		\item The semi-closed channel $[0,2\pi)\times \re$. 
		\item The closed channel $[0,2\pi]\times \re$.
	\end{enumerate}
	
	Let $\rie$ be one of the Riemann surfaces listed in the above. 
	We define two vector fields $V^{\mathrm{rot}},V^{\mathrm{tra}}\in\mathfrak{X}(\rie)$ 
	\begin{align}
		V^{\mathrm{rot}}=\partial_\theta,\quad V^{\mathrm{tra}}=\partial_y,
	\end{align} 
	where $z=r\me^{\mi\theta}=x+\mi y\in \rie$. 
	Then, $V^{\mathrm{rot}}$ is a HKVF on $\rie$ if $\rie= \hat\ce$, $\ce$, $\Delta$, $\ce^*$, $\Delta^*$, $\Delta_\rho$, $\overline\Delta$, $\overline\Delta^*$, $\overline\Delta_\rho$, $\tilde\Delta_\rho$. 
	$V^{\mathrm{tra}}$ is a HKVF on $\rie$ if $\rie= \ce$, $(0,\infty)\times \re$, $(0,2\pi)\times \re$, $\se^1\times \re$, $\te_\Lambda$, $[0,\infty)\times\re$, $[0,2\pi)\times \re$, $[0,2\pi]\times \re$. 
\end{example}

It is known that for every surface, every HKVF on the surface is reduced to one of the above cases via a conformal mapping. 	
\begin{fact}[\cite{Shimizu_2023_Killing}]
\label{fact:main}
	Let $(M,g)$ be a surface. 
	Let $X$ be a hydrodynamic Killing vector field on $M$. 
	Then, there exists a Riemann surface $\rie$ and a conformal mapping $\phi:M\to \rie$ with the conformal factor $\lambda:\rie\to (0,\infty)$ 
	such that 
	\begin{enumerate}
		\item if $\mathrm{Per}(X)\ne \emptyset$, $\rie$ is the one of the following Riemann surfaces: 
		$\hat\ce$, 
		$\ce$, 
		$\Delta$, 
		$\ce^*$, 
		$\Delta^*$, 
		$\Delta_\rho$, 
		$\te_\Lambda$, 
		$\overline\Delta$, 
		$\overline\Delta^*$, 
		$\overline\Delta_\rho$, 
		$\tilde\Delta_\rho$ and
		\begin{align}
			(\phi_*X)_t(z)
			&=
			\begin{cases}
				z+\mi t,&\quad \text{if }\rie=\te_\Lambda,\\
				\me^{\mi t}z,&\quad\text{otherwise},
			\end{cases}\\
			(\phi^{-1})^*g
			&=
			\begin{cases}
				\lambda^2(\mathrm{Re}(z))|\md z|^2,&\quad \text{if }\rie=\te_\Lambda,\\
				\lambda^2(|z|)|\md z|^2,&\quad\text{otherwise}.
			\end{cases}
		\end{align}
		\item otherwise, $\rie$ is the one of the following Riemann surfaces: 
		$\ce$, $(0,\infty)\times \re$, $(0,2\pi)\times \re$, $\se^1\times \re$, $\te_\Lambda$, $[0,\infty)\times\re$, $[0,2\pi)\times \re$, $[0,2\pi]\times \re$
		and 
		\begin{align}
			(\phi_*X)_t(z)
			&=z+\mi t,\\
			(\phi^{-1})^*g
			&=\lambda^2(\mathrm{Re}(z) )|\md z|^2.
		\end{align}
	\end{enumerate}
\end{fact}

In particular, by restricting $M$ to $M\setminus\mathrm{Fix}(X)$, we can utilize more simplified description as follows. 
\begin{fact}[\cite{Shimizu_2023_Killing}]
\label{fact-E}
	Let $(M,g)$ be a surface with a hydrodynamic Killing vector field $X\in\mathfrak{X}(M)$. 
	Then, there exists a conformal mapping $\phi:M\setminus\mathrm{Fix}(X)\to E$ such that 
	\begin{enumerate}
		\item $E$ is either $B\times \se^1$, $B\times \re$ or $\te_\Lambda$ for some $1$-manifold $B$.
		\item $g$ is represented as $\lambda^2(x^1)((\md x^1)^2+(\md x^2)^2)$ on $E$.
		\item $X$ is represented as $\partial_2$ on $E$. 
	\end{enumerate} 
	In particular, $\lambda=|X|\circ\phi^{-1}$. 
\end{fact}
In other words, for any surface $(M,g)$ and any HKVF $X\in\mathfrak{X}(M)$, taking $\phi:M\setminus\mathrm{Fix}(X)\to E$, without loss of generality, we can assume that on $E$, $g$ and $X$ are given as 
\begin{align}
	g&=\lambda^2(x)(\md x^2+\md y^2),\\
	X&= \frac{\partial}{\partial y}.
\end{align}
for some smooth function $\lambda:E\to (0,\infty)$. 

\section{Main result}
\label{sec-main} 
In what follows, we construct an analytic formula of a HGF on a surface with a HKVF by using Fact~\ref{fact:main}. 
Since all surfaces in a conformal class share the same HGF, the surfaces on which we have to construct a HGF are the 14 Riemann surfaces listed in Example~\ref{ex-surf}. 
Unless $M$ is a closed surface, our main task is to choose harmonic part of HGF appropriately to satisfy the boundary and end condition. 
For closed surfaces, once dividing HGF into a harmonic part $-\triangle \Phi(x,x_0)=\delta(x;x_0)$ and a metric potential part $\triangle V=1$, we confirm that the function $G(x,x_0)=\Phi(x,x_0)+V(x)/|M|+V(x_0)/|M|$ is smoothly defined on the whole space except for $x_0$. 
It is hard to find analytic formulae of a HGF but once we found, it is straightforward to check the formulae become a HGF on the surface.
For non-closed surfaces, the following is merely a rewrite of the classically known results by the local coordinate, but for closed surfaces, our formula is worth an extension of one that has been already obtained in some special cases. 

Given a surface $(M,g)$ with a HKVF $X\in\vf(M)$, take the conformal mapping $\phi:M\to \rie$ given in Fact~\ref{fact:main}. 
For simplicity of notation, we use the same notation for a function and its local coordinate representation. 
Let 
\begin{align}
	P(z)=(1-z)\prod_{n\ge 1}(1-\rho^nz )(1-\rho^n z^{-1}). 
\end{align}
The function $P$ is known as a kind of Schottky-Klein prime function and has the following notable property. 
\begin{align}
\label{eq-autoP}
	\log|P(\rho^k z)|=-k\log|z|-\frac{k(k-1)}{2}\log|\rho|+\log|P(z)| 
\end{align}
for each $k\in\ze$.
Set an associated constant $c_P$ by 
\begin{align}
	c_P
	\coloneqq \lim_{z\to 1}\frac{P(z)}{1-z}
	=\prod_{n\ge 1}(1-\rho^n)^2. 
\end{align}
The metric potential of our analytic formula is given as follows."
\begin{align}
	V(z)= \int_1^{|z|}\frac{1}{u} \int_1^u s\lambda^2(s)\md s\md u. 
\end{align}
for each $z\in \rie$. 
$V$ now satisfies 
\begin{align}
\label{eq-V}
	\triangle V=1
\end{align}
since 
\begin{align}
	\triangle V= \lambda^{-2}\left(\frac{\partial^2}{\partial r^2}+\frac{1}{r}\frac{\partial}{\partial r}+\frac{1}{r^2}\frac{\partial^2}{\partial\theta^2} \right)
	\int_1^{r}\frac{1}{u} \int_1^u s\lambda^2(s)\md s\md u=1 
\end{align}
with $z=r\me^{\mi\theta}$. 
Choose a constant $c_e\in\re $ arbitrarily. 

\begin{thm}
\label{thm-Analytic-form}
Let $(M,g)$ be a surface with a hydrodynamic Killing vector field $X\in\vf(M)$. 
Let $\phi:M\to \rie$ be a conformal mapping given in Fact~\ref{fact:main}. 
Then, there exists a hydrodynamic Green function $G_M:M\times M\setminus \{x=x_0\}\to \re$ and a Robin function $R_M:M\to \re$ such that the following functions $G:\rie\times \rie\setminus  \{z=z_0\}\to \re $ and $R:\rie\to \re $ defined on $\rie$ is the local representation of $G_M$ and $R_M$ by $\phi$, respectively.  
\begin{enumerate}
	\item [$\hat \ce$.]
	\begin{align}
		G(z,z_0) 
		=&-\frac{1}{2\pi}\left( \log|z-z_0|-c_v\log|z z_0|\right)\\
		&+\frac{1}{|M|}\int_1^{|z|}\frac{1}{u} \int_1^u s\lambda^2(s)\md s\md u
		+\frac{1}{|M|}\int_1^{|z_0|}\frac{1}{u} \int_1^u s\lambda^2(s)\md s\md u,\\
		R(z)
		=&\frac{1}{2\pi}\left( \log\lambda(z) +c_v\log|z|^2\right)
		+\frac{2}{|M|}\int_1^{|z|}\frac{1}{u} \int_1^u s\lambda^2(s)\md s\md u,
	\end{align} 
	where $c_v\in(0,1)$ is a constant dependent only on the choice of the coordinate, that is, 
	\begin{align}
	c_v =\frac{2\pi }{|M|} \int_0^1 s\lambda^2(s)\md s=1-\frac{2\pi }{|M|} \int_1^\infty s\lambda^2(s)\md s.
	\end{align}
	$c_v$ indicates the rate of the volume for $\phi^{-1}(\Delta)$ relative to $M$. 
	\item [$\ce$.]
	\begin{align}
		G(z,z_0) 
		=&-\frac{1}{2\pi}\log|z-z_0|,\\
		R(z) 
		=&\frac{1}{2\pi}\log\lambda(z).
	\end{align}
	\item [$\Delta$.]
	\begin{align}
		G(z,z_0) 
		=&-\frac{1}{2\pi}\left( \log|z-z_0|-\log|1-z \bar z_0|\right),\\
		R(z) 
		=&\frac{1}{2\pi}\left( \log\lambda(z) +\log(1-|z|^2)\right).
	\end{align}
	\item [$\ce^*$.]
	\begin{align}
		G(z,z_0)
		=&-\frac{1}{2\pi}\left( \log|z-z_0|- c_e\log|z z_0|\right),\\
		R(z)
		=&\frac{1}{2\pi}\left( \log\lambda(z) +\log|z|^2\right).
	\end{align}
	\item [$\Delta^*$.]
	\begin{align}
		G(z,z_0) 
		=&-\frac{1}{2\pi}\left( \log|z-z_0|-\log|1-z \bar z_0|-c_e \log|z z_0|\right),\\
		R(z) 
		=&\frac{1}{2\pi}\left( \log\lambda(z)+\log(1-|z|^2)+c_e \log|z|^2\right).
	\end{align}
	\item [$\Delta_\rho$.]
	\begin{align}
		G(z,z_0)
		=&-\frac{1}{2\pi}\left(\log\left|z_0P(z/z_0)\right|-\log\left|P(z\bar z_0)\right|\right),\\
		R(z)
		=&\frac{1}{2\pi}\left(\log\lambda(z) -\log c_P+\log P(|z|^2)\right).
	\end{align}
	\item [$\te_\Lambda$.]
	\begin{align}
		G(z,z_0)
		=&-\frac{1}{2\pi}\left(\log\left|\me^{z_0}P(\me^{z-z_0})\right|+\frac{\rp(z)\rp(z_0)}{\rp(\tau)}-\frac{1}{2}(\rp(z)+\rp(z_0))\right)\\
		&+\frac{1}{|M|}\int_0^{\rp(z)} \int_0^u \lambda^2(s)\md s\md u
		+\frac{1}{|M|}\int_0^{\rp(z_0)} \int_0^u \lambda^2(s)\md s\md u,\\
		R(z)
		=&\frac{1}{2\pi}\left(\log\lambda(z) -\log c_P-\frac{\rp(z)^2}{\rp(\tau)}+\rp(z)\right)
		+\frac{2}{|M|}\int_0^{\rp(z)} \int_0^u \lambda^2(s)\md s\md u,
	\end{align}
	where $\Lambda =\{2\pi\mi m +\tau n|\,\mathrm{Re}(\tau)>0,\,m,n\in\ze \}$
	and $\rho=\exp(-\rp(\tau))\in(0,1)$. 
	\item [$\overline\Delta$.]
	\begin{align}
		G(z,z_0)
		=&-\frac{1}{2\pi}\left( \log|z-z_0|-\log|1-z \bar z_0|\right),\\
		R(z) 
		=&\frac{1}{2\pi}\left( \log\lambda(z) +\log(1-|z|^2)\right).
	\end{align}
	\item [$\overline\Delta^*$.]
	\begin{align}
		G(z,z_0)
		=&-\frac{1}{2\pi}\left( \log|z-z_0|-\log|1-z \bar z_0|-c_e \log|z z_0|\right),\\
		R(z) 
		=&\frac{1}{2\pi}\left( \log\lambda(z)+\log(1-|z|^2)+c_e \log|z|^2\right).
	\end{align}
	\item [$\overline\Delta_\rho$.]
	\begin{align}
		G(z,z_0)
		=&-\frac{1}{2\pi}\left(\log\left|z_0P(z/z_0)\right|-\log\left|P(z\bar z_0)\right|\right),\\
		R(z)
		=&\frac{1}{2\pi}\left(\log\lambda(z) -\log c_P+\log P(|z|^2)\right).
	\end{align}
	\item [$\tilde\Delta_\rho$.]
	\begin{align}
		G(z,z_0)
		=&-\frac{1}{2\pi}\left(\log\left|z_0P(z/z_0)\right|-\log\left|P(z\bar z_0)\right|\right),\\
		R(z)
		=&\frac{1}{2\pi}\left(\log\lambda(z) -\log c_P+\log P(|z|^2)\right).
	\end{align}
	\item [$[0,\infty)\times\re$.]
	\begin{align}
		G(z,z_0)
		=&-\frac{1}{2\pi}\left( \log|z-z_0|-\log|z+\bar z_0|\right),\\
		R(z)
		=&\frac{1}{2\pi}\left(\log\lambda(z) +\log |z+\bar z|\right).
	\end{align}	
	\item [$[0,2\pi)\times\re$.]
	\begin{align}
		G(z,z_0)
		=&-\frac{1}{2\pi}\left( \log\left|\me^{\mi z/2} -\me^{\mi z_0/2}\right|-\log\left|\me^{\mi z/2} -\me^{-\mi \bar z_0/2}\right|
		-c_e(\rp(z)+\rp(z_0))
		\right),\\
		R(z)
		=&\frac{1}{2\pi}\left( \log\lambda(z)-\frac{1}{2}\log\left|\me^{\mi z/2}\right|+\log\left|\me^{\mi z/2} -\me^{-\mi \bar z/2}\right|
		+2c_e\rp(z)
		\right). 
	\end{align}
	\item [${[0,2\pi]\times\re}$.]
	\begin{align}
		G(z,z_0)
		=&-\frac{1}{2\pi}\left( \log\left|\me^{\mi z/2} -\me^{\mi z_0/2}\right|-\log\left|\me^{\mi z/2} -\me^{-\mi \bar z_0/2}\right|
		-c_e(\rp(z)+\rp(z_0))
		\right),\\
		R(z)
		=&\frac{1}{2\pi}\left( \log\lambda(z)-\frac{1}{2}\log\left|\me^{\mi z/2}\right|+\log\left|\me^{\mi z/2} -\me^{-\mi \bar z/2}\right|
		+2c_e\rp(z)
		\right). 
	\end{align}
\end{enumerate}
\end{thm}

\begin{proof}
The proof differs depending on whether the surface is a closed surface or not. 
If the surface $(M,g)$ is a closed surface, we confirm that $G$ is an HGF on $(\rie,\lambda^2|\md z|^2)$. 
Otherwise, owing to the conformal invariance of HGFs, it is sufficient to see that $G$ is an HGF on $(\rie,|\md z|^2)$. 
The majority of the proof is devoted to demonstrating that when $\rie$ is a closed surface, that is, $\rie=\hat\ce,\te_\Lambda$, $G$ can be smoothly defined on $\rie \times \rie \setminus\{z=z_0\}$. 
For the rest, what needs to be shown can be straightforwardly confirmed through the following procedure.

First, let us see that the case where $\rie$ is not a closed surface, $G$ is smooth on $\rie \times \rie \setminus \{z = z_0\}$ and satisfies Eq.~\eqref{eq-para-end} and Eq.~\eqref{eq-hyp-end}. 
In any relevant case, $G$ can be continuously extended to $\{|z|=1\}$ and consequently, Eq.~\eqref{eq-hyp-end} is satisfied. 
It is clear for $\rie=\ce$. 
For $\rie=\Delta,\overline\Delta$, this follows from 
\begin{align}
\label{eq-disc-bd}
	-\frac{1}{2\pi}\left( \log|z-z_0|-\log|1-z \bar z_0|\right)
	=0	
\end{align}
for each $|z|=1$. 
For $\rie=\ce^*,\Delta^*,\overline\Delta^*$, since the limit 
\begin{align}
	\lim_{|z|\to 0}\left( G(z,z_0)-\frac{c_e}{2\pi }\log|z| \right)
\end{align}
exists, Eq.~\eqref{eq-para-end} is satisfied. 
For $\rie=\Delta_\rho,\overline\Delta_\rho,\tilde\Delta_\rho$, we can deduce from Eq.~\eqref{eq-autoP} 
\begin{align}
	G(\rho z,z_0)= G(z,z_0)-\pi^{-1}\log|z_0|,
\end{align}
which implies 
\begin{align}
	G(z,z_0)=0,\quad 
	G(\rho z,z_0)=-\pi^{-1}\log|z_0|
\end{align} 
for each $|z|=1$. 
For $\rie=[0,\infty)\times\re$, if $\rp(z)=0$, it is obvious that  
\begin{align}
	G(z,z_0)=0
\end{align}
owing to $\bar z=-z$. 
For $[0,2\pi)\times\re,[0,2\pi]\times\re$, if $\rp(z)\in \{0,2\pi\}$, it follows from 
\begin{align}
	\log\left|\me^{\mi z/2} -\me^{\mi z_0/2}\right|=\log\left|\me^{\mi z/2} -\me^{-\mi \bar z_0/2}\right|
\end{align}
that 
\begin{align}
	G(z,z_0)=
	\begin{cases}
		c_e\rp(z_0)/2\pi ,\quad &\text{if }\rp(z)=0,\\
		c_e(2\pi +\rp(z_0))/2\pi, \quad &\text{if }\rp(z)=2\pi. 
	\end{cases}
\end{align}
Consequently, $G$ is smooth on $\rie \times \rie \setminus \{z = z_0\}$. In all cases, the differentiability of $G$ is derived from the representation $G = \log |h|$, where $h$ is a holomorphic function $h: \rie \setminus \{z_0\} \to \ce$.

Once this is shown, for proving that $G$ is an HGF, it remains only to show that $G$ is the Green function for the Laplacian $\triangle=4\lambda^{-2}\partial\bar\partial$, which is straightforward. 
Since the harmonic part of our analytic formulae is assembled by the following functions, 
\begin{align}
\label{eq-harmonics}
	-\triangle \log|z-z_0|
	&=2\pi\delta(z,z_0),\\
	-\triangle \log|z|
	&=2\pi\delta(z,0),\\
	-\triangle \log|1-z\bar z_0|
	&=2\pi\delta(z,\bar z_0^{-1}),\\
	-\triangle \log|z+\bar z_0|
	&=2\pi\delta(z,-\bar z_0),
\end{align}
except for $R=\Delta_\rho$, $\te_\Lambda$, $\overline\Delta_\rho$, $\tilde\Delta_\rho$, $[0,2\pi)\times\re$, $[0,2\pi]\times\re$, we then see that $G$ is an HGF. 
However, in the case of $R=\hat \ce$, it remains unconfirmed whether $-\triangle G(z,z_0)=\delta(z,z_0)-|M|^{-1}$ also holds for $z=0,\infty$, which will be addressed later on. 

For $R=\Delta_\rho$, $\overline\Delta_\rho$ and $\tilde\Delta_\rho$, since 
\begin{align}
	 \log| z_0 P(z/z_0)|
	&=\log|z-z_0|+\sum_{n\ge1}\left(\log\left|1-\rho^n\frac{z}{z_0}\right|+\log\left|1-\rho^n\frac{z_0}{z}\right| \right),\\
	\log| P(z\bar z_0)|
	&=\log|1-z\bar z_0|+\sum_{n\ge1}\left(\log\left|1-\rho^n z\bar z_0\right|+\log\left|1-\rho^nz^{-1}\bar z_0^{-1}\right| \right),
\end{align}
we obtain
\begin{align}
	-\triangle \log| z_0 P(z/z_0)| &=-\triangle \log|z-z_0|=2\pi\delta(z,z_0)\\
	-\triangle \log |P(z\bar z_0)|&= 0,
\end{align}
which yields $G$ is an HGF.
By choosing $\tilde\Delta_\rho$ as a fundamental domain for $\te_\Lambda$, we can deduce that $G$ is an HGF on $\te_\Lambda$. 
For the remaining two cases, $[0,2\pi)\times\re$ and $[0,2\pi]\times\re$, it is confirmed from 
\begin{align}
	-\triangle\left(\log\left|\me^{\mi z/2} -\me^{-\mi \bar z_0/2}\right|
		-c_e(\rp(z)+\rp(z_0))\right)=0
\end{align}
and 
\begin{align}
	 -\triangle \log\left|\me^{\mi z/2} -\me^{\mi z_0/2}\right|
	 &=-\triangle \log |z-z_0|-\triangle \log\left|\frac{\me^{\mi z/2} -\me^{\mi z_0/2}}{z-z_0}\right|\\
	 &=2\pi \delta(z,z_0),
\end{align}
since $f(z)=(\me^{\mi z/2} -\me^{\mi z_0/2})/(z-z_0)$ is a holomorphic function with removable singularity at $z=z_0$. 

After proving that $G$ is an HGF, $R$ can be straightforwardly derived from
\begin{align}
\label{eq-R-G}
	R(z)=\lim_{z_0\to z} \left(G(z,z_0)+\frac{1}{2\pi}\log |z-z_0|\right)+\frac{1}{2\pi} \log \lambda(z)
\end{align}
owing to $d(z,z_0)=\lambda(z)|z-z_0|+o(1)$. 
The remaining task is to show that, in the case of $\rie = \hat{\ce}$ or $\te_\Lambda$, $G$ is smoothly defined on $\rie \times \rie \setminus \{z = z_0\}$, and in particular, for $\rie = \hat{\ce}$, $-\triangle G(z, z_0) = \delta(z, z_0) - |M|^{-1}$ also holds for $z = 0, \infty$.

\begin{enumerate}
	\item [$\hat \ce$.]
	We first see $G$ is smoothly defined, especially, at $z=0$ and $\infty$. 
	Concerning $z=0$, we employ the following identity. 
	\begin{align}
	\label{eq-z=0}
		\frac{1}{|M|}\int_1^{|z|}\frac{1}{u} \int_1^u s\lambda^2(s)\md s\md u= \frac{1}{|M|}\int_1^{|z|}\frac{1}{u} \int_0^u s\lambda^2(s)\md s\md u-\frac{c_v}{2\pi}\log|z|
	\end{align}
	since
	\begin{align}
		 \int_1^u s\lambda^2(s)\md s
		 &=\int_0^u s\lambda^2(s)\md s-\int_0^1 s\lambda^2(s)\md s\\
		 &=\int_0^u s\lambda^2(s)\md s-\frac{c_v|M|}{2\pi}. 
	\end{align}
	It follows from Eq.~\eqref{eq-z=0} that 
	\begin{align}
		G(z,z_0) 
		=&-\frac{1}{2\pi}\left( \log|z-z_0|-c_v\log|z z_0|\right)\\
		&+\frac{1}{|M|}\int_1^{|z|}\frac{1}{u} \int_0^u s\lambda^2(s)\md s\md u-\frac{c_v}{2\pi}\log|z|\\
		&+\frac{1}{|M|}\int_1^{|z_0|}\frac{1}{u} \int_1^u s\lambda^2(s)\md s\md u\\
		=&-\frac{1}{2\pi}\left( \log|z-z_0|-c_v\log|z_0|\right)\\
		&+\frac{1}{|M|}\int_1^{|z|}\frac{1}{u} \int_0^u s\lambda^2(s)\md s\md u
		+\frac{1}{|M|}\int_1^{|z_0|}\frac{1}{u} \int_1^u s\lambda^2(s)\md s\md u. \label{eq-G-z=0}
	\end{align}
	Owing to $\lambda$ of class $C^1$, we have 
	\begin{align}
		\lambda^2(s)=\lambda^2(0)+o(1),
	\end{align}
	which gives 
	\begin{align}
		\int_1^{|z|}\frac{1}{u} \int_0^u s\lambda^2(s)\md s\md u
		=\frac{\lambda^2(0)}{4}(|z|^2-1)+o(|z|^2). 
	\end{align}
	Hence, $G$ is continuous at $z=0$ and $-\triangle G(z,z_0)=\delta(z,z_0)-|M|^{-1}$ since 
	\begin{align}
		\triangle \left(\int_1^{|z|}\frac{1}{u} \int_0^u s\lambda^2(s)\md s\md u\right)=1.
	\end{align}
	In the same way, we deduce the smoothness of $G$ at $z=0$. 
	Concerning $z=\infty$, since
	\begin{align}
		\int_1^u s\lambda^2(s)\md s
		 &=\int_\infty^u s\lambda^2(s)\md s+\int_1^\infty s\lambda^2(s)\md s\\
		 &=\int_\infty^u s\lambda^2(s)\md s-\frac{(c_v-1)|M|}{2\pi}, 
	\end{align}
	we obtain
	\begin{align}
		G(z,z_0) 
		=&-\frac{1}{2\pi}\left( \log|z-z_0|-c_v\log|z z_0|\right)\\
		&+\frac{1}{|M|}\int_1^{|z|}\frac{1}{u} \int_\infty^u s\lambda^2(s)\md s\md u-\frac{c_v-1}{2\pi}\log|z|\\
		&+\frac{1}{|M|}\int_1^{|z_0|}\frac{1}{u} \int_1^u s\lambda^2(s)\md s\md u\\
		=&-\frac{1}{2\pi}\left( \log\left|1-\frac{z_0}{z}\right|-c_v\log|z_0|\right)\\
		&+\frac{1}{|M|}\int_1^{|z|}\frac{1}{u} \int_\infty^u s\lambda^2(s)\md s\md u
		+\frac{1}{|M|}\int_1^{|z_0|}\frac{1}{u} \int_1^u s\lambda^2(s)\md s\md u.\label{eq-G-z=infty}
	\end{align}
	Hence, Using the same argument as in the case at $z=0$, the smoothness of $G$ at $z=\infty$ is obtained. 
	
	It follows from Eq.~\eqref{eq-harmonics} and \eqref{eq-V} that $G$ satisfies $-\triangle G(z,z_0)=\delta(z,z_0) -|M|^{-1}$, except at $z=0,\infty$. 
	For $z=0,\infty$, this is shown by using the equivalent presentations for $G$, given by \eqref{eq-G-z=0} and \eqref{eq-G-z=infty}, which yields $G$ is a HGF. 
	Since 
	\begin{align}
		\lim_{z_0\to z}\left(G(z,z_0)+\frac{1}{2\pi}\log |z-z_0|\right)
		&=\frac{c_v}{2\pi}\log|z|^2+\frac{2}{|M|}\int_1^{|z|}\frac{1}{u} \int_1^u s\lambda^2(s)\md s\md u,
	\end{align}
	using Eq.~\eqref{eq-R-G}, we obtain the analytic formula of $R$. 
	\item [$\te_\Lambda$.] 
	To demonstrate the smoothness of $G$, it suffices to show that for each $\gamma, \gamma_0 \in \Lambda$,
    \begin{align}
        G(z + \gamma, z + \gamma_0) = G(z, z_0).
    \end{align}
    Set 
    \begin{align}
    	\gamma &=2\pi \mi m+\tau n,\quad \gamma_0=2\pi \mi m_0+\tau n_0,\\
    	x&=\rp(z),\quad x_0=\rp(z_0),\quad T=\rp(\tau). 
    \end{align}
	Owing to 
	\begin{align}
		\me^{z+\gamma-z_0-\gamma_0}=\me^{(n-n_0)\rp(\tau)+z-z_0}=\rho^{n_0-n}\me^{z-z_0},
	\end{align}
	we deduce from Eq.~\eqref{eq-autoP} that 
	\begin{align}
		&\log|\me^{z_0+\gamma_0} P(\rho^{n_0-n}\me^{z-z_0})|-\log |\me^{z_0}P(\me^{z-z_0})|\\
		=&n_0T-(n_0-n)(x-x_0)+\frac{1}{2}(n_0-n)(n_0-n-1)T.
	\end{align}
	Since $\lambda^2$ is a $T$-periodic function, the Fourier expansion 
	\begin{align}
		\lambda^2(s)=a_0+\sum_{n\ge1}a_n\cos(\omega_n s)+b_n\sin(\omega_n s),
	\end{align}
	with $\omega_n=2\pi n /T$ yields that 
	\begin{align}
		\int_0^{x}\int_0^u\lambda^2(s)\md s\md u
		=\frac{1}{2}a_0 x^2-\sum_{n\ge1}\omega_n^{-2}(a_n\cos(\omega_n x)+b_n\sin(\omega_n x))
	\end{align}
	Note that 
	\begin{align}
		|M|=\int_0^{2\pi}\int_0^T\lambda^2(s)\md s\md u,
	\end{align}
	which gives
	\begin{align}
		\int_0^T\lambda^2(s)\md s=a_0 T=\frac{|M|}{2\pi }.
	\end{align}
	Hencem, we obtain  
	\begin{align}
		\frac{1}{|M|}\left(\int_0^{\rp(z+\gamma)}\int_0^u\lambda^2(s)\md s\md u-\int_0^{\rp(z)}\int_0^u\lambda^2(s)\md s\md u\right)
		=\frac{1}{2\pi}n\left(x+\frac{nT}{2}\right).
	\end{align}
	Since
	\begin{align}
		\frac{\rp(z+\gamma)\rp(z_0+\gamma_0)}{\rp(\tau)}-\frac{\rp(z)\rp(z_0)}{\rp(\tau)}
		&=n_0 x+n x_0+nn_0 T,
	\end{align}
	we obtain 
	\begin{align}
		&2\pi(G(z+\gamma,z_0+\gamma_0)-G(z,z_0))\\
		=&-n_0T+(n_0-n)(x-x_0)-\frac{1}{2}(n_0-n)(n_0-n-1)T\\
		&-n_0 x-n x_0-nn_0 T+\frac{1}{2}(n+n_0)T\\
		&+n\left(x+\frac{nT}{2}\right)+n_0\left(x_0+\frac{n_0T}{2}\right)\\
		=&0.
	\end{align}
\end{enumerate}

\end{proof}

\section{Impact of the curvature on a point vortex}
\label{sec-impact-curvature}
\subsection{General relation between curvature and point vortex dynamics}
\label{subsec-general-curvature}
As an application of Theorem\ref{thm-Analytic-form}, we clarify the mechanism whereby the curvature affects a point vortex in both qualitative and quantitative viewpoints. 
To begin with the qualitative aspect, we prove the following relation between the Robin function $R$ and the Gauss curvature $\kappa$.
\begin{prop}
Let $(M,g)$ be a surface. Then, 
	\begin{align}
	\label{eq-Robin-curvature}
		-\triangle R= \frac{\kappa}{2\pi}+2c_M.
	\end{align}
\end{prop}
\begin{proof}
	Let $z$ is a complex coordinate with a conformal factor $\lambda$ centered at $x_0\in M$ with $z=z(x)$, $z_0=z(x_0)$. 
	Take a hydrodynamic Green function $G\in C^\infty(M\times M\setminus\Delta)$. 
	In the coordinate $z$, $G$ is presented as 
	\begin{align}
		G(z,z_0)=-\frac{1}{2\pi}\log |z-z_0|+h(z)+h(z_0)+V(z)+V(z_0)
	\end{align}
	for some harmonic function $h$ and metric potential $V$
	\begin{align}
		-\triangle h&=0,\\
		-\triangle V&=c_M.
	\end{align}
	Since the distance $d$ satisfies 
	\begin{align}
		d(z,z_0)=\lambda(z)|z-z_0|+o(|z-z_0|),
	\end{align}
	we deduce 
	\begin{align}
		R(z)&=\lim_{z_0\to z}\left( G(z,z_0)+\frac{1}{2\pi}\log d(z,z_0)\right)\\
		&= \frac{1}{2\pi}\log \lambda(z)+2h(z)+2V(z)
	\end{align}
	Owing to 
	\begin{align}
		\kappa=-\triangle \log \lambda,
	\end{align}
	we obtain 
	\begin{align}
		-\triangle R= \frac{\kappa}{2\pi}+2c_M.
	\end{align}
\end{proof}
From Eq.~\eqref{eq-jgradR}, the motion of a point vortex is described by the Hamiltonian vector field $X_H = -\jgrad R$, with the Robin function $R$ serving as the Hamiltonian.
In particular, when $M$ is not a closed surface, Eq.~\eqref{eq-Robin-curvature} shows that the vorticity is directly given by the curvature.
Consequently, the point vortex dynamics can be understood precisely as a Hamiltonian flow whose vorticity is dictated by the Gaussian curvature. Explicitly, the curvature drives the point vortex $q$ according to:
\begin{align}
\label{eq-1-vortex-curvature}
\dot q=-\jgrad (-\triangle)^{-1}\left( \frac{\kappa}{2\pi}+2c_M\right).
\end{align}
At first glance, Eq.~\eqref{eq-1-vortex-curvature} appears as a linear response mediated by the inverse Laplacian operator $(-\triangle)^{-1}$.
However, the Riemannian metric corresponding to a prescribed Gaussian curvature is defined as a solution of a semilinear elliptic PDE known as the prescribed curvature problem (or Berge problem).
Therefore, the influence of curvature on the point vortex inherently involves nonlinear effects.
Thus, although Eq.~\eqref{eq-1-vortex-curvature} succinctly and explicitly describes the relationship between curvature and vortex motion, it is generally not easy to analyze in detail.
Nevertheless, in the present work, we succeeded in deriving analytic solutions explicitly (which constitute the main results of this paper and the basis for the analysis conducted in this section), thereby enabling a detailed mathematical characterization of various scaling limits and clearer understanding of the underlying mathematical structure.
Additionally, when $M$ is a closed surface, one must account not only for curvature but also for the global geometric effects encapsulated by the constant $c_M$, which is related to the total surface area.

\begin{rem}
Given a manifold $M$, the problem of finding a Riemannian metric on $M$ whose scalar curvature is a prescribed function $R$ is called the \emph{prescribed scalar curvature problem}. In particular, when $M$ is two-dimensional, restricting the search to the conformal class $\tilde g = \mathrm{e}^{2\sigma}g$ reduces the problem to finding a Riemannian metric on the given surface $(M,g)$ whose Gauss curvature $K$---satisfying $K = R/2$ in two dimensions---matches a prescribed function on $M$. 

Then, in order to derive the equation that characterizes $\tilde g = \mathrm{e}^{2\sigma}g$ with Gauss curvature $\tilde K$, we proceed as follows.
Using local isothermal coordinates $(x^1, x^2)$, one may write the metric in the form
\begin{align}
g &= \mathrm{e}^{2f}\left((\md x^1)^2 + (\md x^2)^2\right)
\end{align}
for a smooth function $f$, in which case the Gauss curvature satisfies
\begin{align}
K &= -\mathrm{e}^{-2f}\triangle f.
\end{align}

We now consider a new metric
\begin{align}
\tilde g = \mathrm{e}^{2\sigma}g 
         = \mathrm{e}^{2\left(f + \sigma\right)}\left((\md x^1)^2 + (\md x^2)^2\right).
\end{align}
Its Gauss curvature $\tilde K$ is given by
\begin{align}
\tilde K &= -\mathrm{e}^{-2\left(f + \sigma\right)}\triangle\left(f + \sigma\right).
\end{align}
Then, it can be rewritten as
\begin{align}
-\triangle_g\,\sigma + K & = \tilde K\,\mathrm{e}^{2\sigma},
\end{align}
which is called \emph{Kazdan--Warner equation}.
Hence, in the prescribed Gaussian curvature problem, one seeks to solve 
the above PDE for $\sigma$, given the target function $\tilde K$.

In this way, by solving the Kazdan--Warner equation, one obtains solutions 
to the prescribed Gaussian curvature problem 
(see, for example, \cite{Aubin_1998,Kazdan_Warner_1974}). 
Meanwhile, since the Kazdan--Warner equation is a semilinear elliptic equation 
with a Liouville-type nonlinearity, one can interpret the influence of curvature 
on the Riemannian metric as a \emph{nonlinear response}.
\end{rem}

\subsection{Analytic formula of the vortex velocity for periodic curvature}
\label{subsec-periodic-curvature}
Building on the qualitative observations from the previous section, we now quantitatively examine the effect of curvature on the dynamics of a point vortex in greater detail. 
To investigate the role of the area constant explicitly, we must focus on closed surfaces.
Among the possible closed surfaces, the analytic formulas established by Theorem~\ref{thm-Analytic-form} are applicable explicitly to either a sphere or a torus; here, we choose a torus, which possesses a richer and more nontrivial topology, as the underlying Riemann surface.
For simplicity, we assume the lattice defining this torus to be rectangular: 
$M=\te_\Lambda$ with $\Lambda =\{2\pi \mi m+ Tn\mid m,n\in\ze \}$ for some $T\in (0,\infty)$.

Given that the curvature on a torus is inherently periodic, when considering its Fourier transform, it is reasonable to regard the scenario in which the curvature is expressed by a one-variable function consisting only of a constant term and the dominant mode as the most fundamental yet nontrivial case.
\begin{align}
	\kappa(x,y)=\kappa_0+\kappa_1\cos \omega x
\end{align}
When $\kappa_0$, $\kappa_1$ and $\omega$ satisfy
$\kappa_0\in (-\infty,0)$, $\kappa_1\in(-\kappa_0,\infty)$, $\omega=2\pi N/T$ for some $N\in\ze_{>0}$, while it is generally difficult to obtain analytic solutions of the Kazdan–Warner equation due to its nonlinearity, we fortunately obtain an explicit analytic solution for the metric $g=\lambda^2|\md z|^2$ that realizes this curvature distribution.
\begin{align}
	\lambda=\frac{\omega\sqrt{-\kappa_0}}{\kappa_1+\kappa_0\cos \omega x}
\end{align}
This result can be readily confirmed through the direct computation below.
\begin{align}
	-\triangle \log \lambda
	&=-\lambda^{-2}(\log\lambda)''\\
	&=-\frac{(\kappa_1+\kappa_0\cos \omega x)^2}{\omega^2(-\kappa_0)}\frac{\omega^2\kappa_0(\kappa_0+\kappa_1\cos \omega x)}{(\kappa_1+\kappa_0\cos \omega x)^2}\\
	&=\kappa_0+\kappa_1\cos \omega x=\kappa.
\end{align}

Define $E:\re\to \re$ by
\begin{align}
E(x)=\int_0^x\lambda^2(s)\md s.
\end{align}
Then, straightforward integration gives
\begin{align}
E(x)=\frac{1}{\kappa_1^2-\kappa_0^2}
\left(\frac{\omega\kappa_0^2\sin \omega x}{\kappa_1+\kappa_0\cos \omega x}-\omega^2\kappa_0\kappa_1 \int_0^x \frac{\md s}{\kappa_1+\kappa_0\cos \omega s}\right),
\end{align}
with the integral expressed explicitly as
\begin{align}
	\int_0^x \frac{1}{\kappa_1+\kappa_0\cos \omega s}\md s
	&=
	\begin{dcases}
		\frac{2}{\omega\sqrt{\kappa_1^2-\kappa_0^2}} 
		\left(\frac{n\pi}{2}+\tan^{-1}\left(\sqrt{\frac{\kappa_1-\kappa_0}{\kappa_1+\kappa_0}}\tan\left(\frac{\omega x}{2}\right)\right)\right)
		,\quad &\text{if $n$ is even},\\
		\frac{2}{\omega\sqrt{\kappa_1^2-\kappa_0^2}} 
		\left(\frac{n\pi}{2}-\tan^{-1}\left(\sqrt{\frac{\kappa_1+\kappa_0}{\kappa_1-\kappa_0}}\cot\left(\frac{\omega x}{2}\right)\right)\right)
		,\quad &\text{if $n$ is odd}.
	\end{dcases}
\end{align}
where $n=\lfloor \omega x/\pi\rfloor$, $x_0=x-n\pi/\omega$.
Consequently, one has
\begin{align}
E(T)&=\frac{T\omega^2 (-\kappa_0) \kappa_1 }{(\kappa_1^2-\kappa_0^2)^{3/2}},\\
|M|&=\int_0^{2\pi}\int_0^T\lambda^2(x)\md x\md y=2\pi E(T).
\end{align} 

Applying Theorem~\ref{thm-Analytic-form}, the Robin function on $M$ explicitly reads
\begin{align}
R(x)=\frac{1}{2\pi}\left(\log\lambda -\log c_P-\frac{x^2}{T}+x\right)+\frac{2}{|M|}\int_0^{x}E(u)\md u.
\end{align}
Since $R$ is independent of $y$, the vortex velocity $\dot{q}_y$ is explicitly obtained by
\begin{align}
\dot q_y
=-\frac{\lambda^{-2}}{2\pi }\left((\log\lambda)'-\frac{2x}{T}+1 +\frac{2E(x)}{E(T)} \right),
\end{align}
with $(\log\lambda)'=\omega \kappa_0\sin\omega x/(\kappa_1+\kappa_0\cos \omega x)$.

\subsection{Linear response for small background curvature}
\label{subsec-linear-response}
\begin{figure}[htbp]
    \centering
    \includegraphics[width=0.3\textwidth]{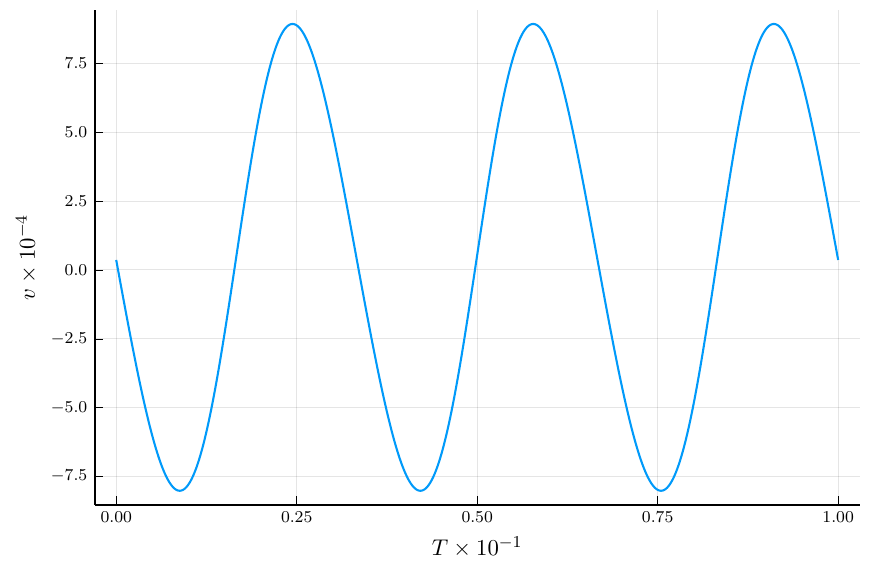}
    \includegraphics[width=0.3\textwidth]{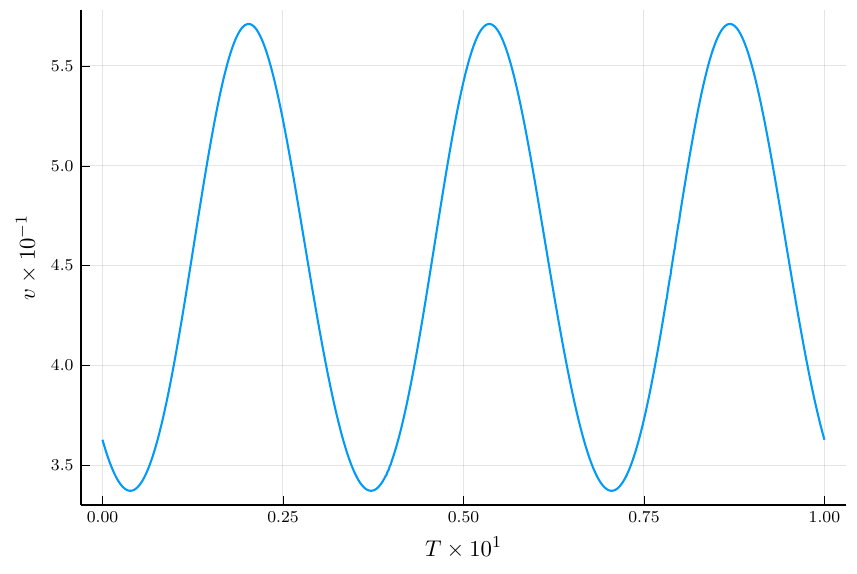}
    \includegraphics[width=0.3\textwidth]{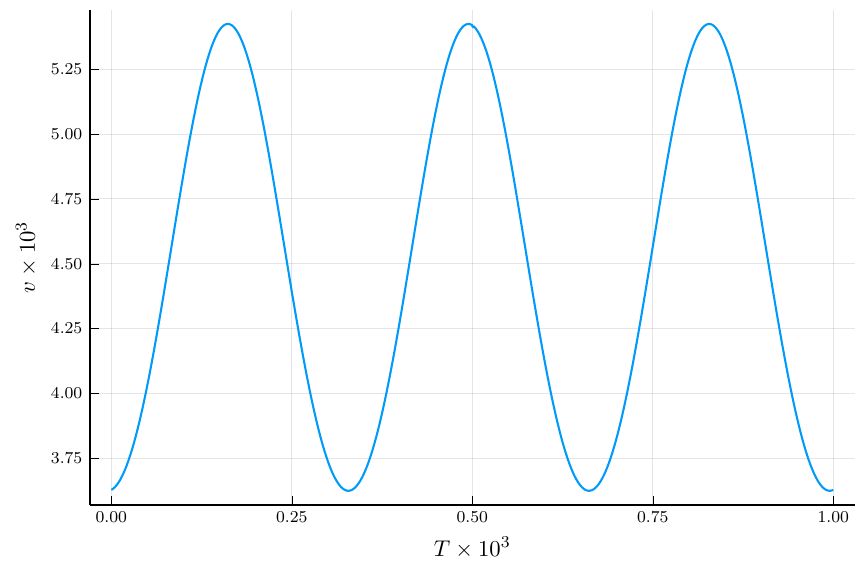}
    \caption{(Left) \(T = 10^{-1}\), (Center) \(T = 10\), (Right) \(T = 10^{3}\).}
    \label{fig-small-kappa0}
\end{figure}

The discussion below focuses on the waveform of $v=\dot{q}_y$ in accordance with the scale of each parameter.
First, let us consider the case where \(\kappa_0\) is sufficiently small.
Figure~\ref{fig-small-kappa0} illustrates the behavior of the system with
\(\kappa_0 = 10^{-1}\), \(\kappa_1 = 1\), and \(\omega = 2\pi \times 3 / T\)
for different time scales \(T\). We compare the results for 
\(T = 10^{-1}\), \(T = 10\), and \(T = 10^{3}\) to highlight how changing
the time scale influences the system's response.
The figure illustrates that, when $\kappa_0 \to 0$, the velocity $v$ exhibits a single waveform having the same period as the curvature while the phase shift is scale-dependent on $T$

To analytically clarify the linear response observed above, we perform a straightforward asymptotic expansion as $\kappa_0 \to 0$. Although the following calculation is elementary, we explicitly show the steps here for completeness and ease of verification.
First, expanding each term individually, we obtain
\begin{align}
\lambda^{-2}
&= -\frac{\kappa_1^2}{\kappa_0\omega^2}\left(1+\frac{2\kappa_0}{\kappa_1}\cos \omega x\right)+o(1), \\
(\log \lambda)'
&= \frac{\omega\kappa_0}{\kappa_1}\sin \omega x+o(\kappa_0).
\end{align}

Next, the integral $E(x)$ can be expanded straightforwardly via a Taylor series as
\begin{align}
E(x)
&= \frac{\omega\kappa_0}{\kappa_1^2}\left(-\omega x+\frac{\omega\kappa_0}{\kappa_1}\sin\omega x\right)+o(\kappa_0^2), \\
E(T)
&= \frac{\omega\kappa_0}{\kappa_1^2}(-\omega T)+o(\kappa_0^2).
\end{align}
Thus, the ratio $2E(x)/E(T)$ simplifies neatly to
\begin{align}
\frac{2E(x)}{E(T)}
&= \frac{2x}{T}-\frac{2\kappa_0}{T\kappa_1}\sin\omega x+o(\kappa_0).
\end{align}

Combining these expansions, we explicitly arrive at the asymptotic expression for the velocity field:
\begin{align}
v = \frac{\kappa_1}{2\pi \omega^2}\left(\frac{\kappa_1}{\kappa_0}+2\cos \omega x+\left(\omega-\frac{4}{T}\right) \sin \omega x\right)+o(1), \quad \text{as } \kappa_0\to 0.
\end{align}

This concise expression clearly demonstrates the linear response, explicitly highlighting how each parameter contributes to the structure and scaling of the velocity waveform observed previously.
A noteworthy point here is that the output velocity $v$ shares the same periodicity as the input curvature oscillation $\cos(\omega x)$. 
Specifically, $v$ can be explicitly expressed as a linear combination of $\sin(\omega x)$, $\cos(\omega x)$, and a constant.
Consequently, it is reasonable to infer that, as $\kappa_0\to 0$, the response of $v$ to the input curvature $\kappa$ is linear.
The primary highlight of the above analysis is the affine similarity law between the waveforms of $v$ and $\kappa$:
the velocity waveform is similar to the curvature waveform, with a similarity ratio  
\begin{align} 
	\frac{1}{\pi \omega^2}\sqrt{1+\left(\frac{\omega}{2}-\frac{2}{T}\right)^2}. 
\end{align}
and phase shift
\begin{align}
	\tan^{-1}\left(\frac{\omega}{2}-\frac{2}{T}\right).
\end{align}

An additional notable point is that
the magnitude of the constant component of the velocity is determined by the representative scale 
\begin{align} 
	\frac{\kappa_1^2}{2\pi \omega^2\kappa_0}. 
\end{align} 
Also, for sufficiently small periods $T$, the phase shift between the velocity and curvature waveforms approaches $+\pi/2$, while becoming smaller as the period $T$ becomes larger.

\subsection{Nonlinear response near resonance of curvature amplitudes}
\label{subsec-nonlinear-response}

Let us take a sufficiently small positive number $\varepsilon$ and set the curvature parameters as
\[
\kappa_1 = K(1+\varepsilon),\quad \kappa_0 = -K\quad (K>0).
\]
This choice implies that the system approaches the critical condition $|\kappa_1| - |\kappa_0| = 0$, corresponding to a boundary situation where the fundamental assumption $|\kappa_1| - |\kappa_0| > 0$ for model validity is nearly violated. Thus, the nonlinearity explicitly emerges in this regime.
\begin{figure}[htbp]
    \centering
    \includegraphics[width=0.3\textwidth]{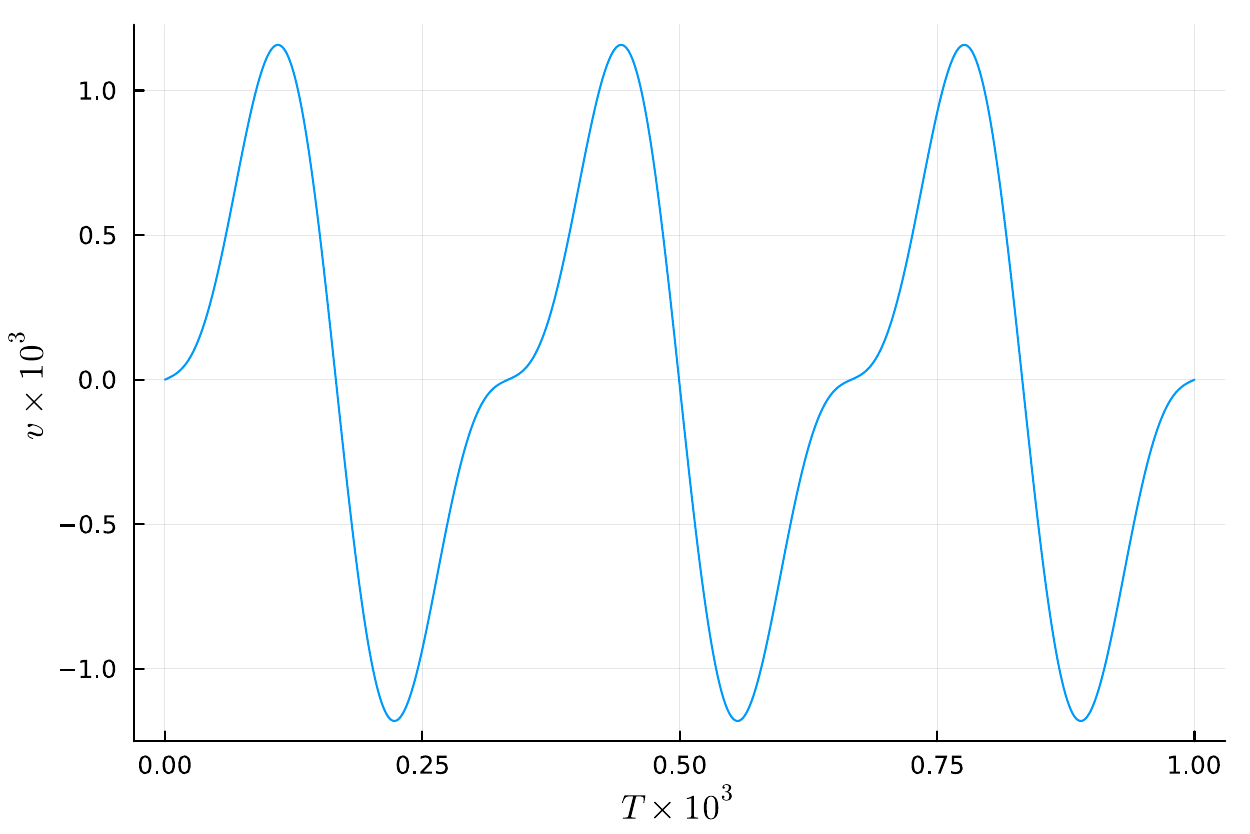}
    \includegraphics[width=0.3\textwidth]{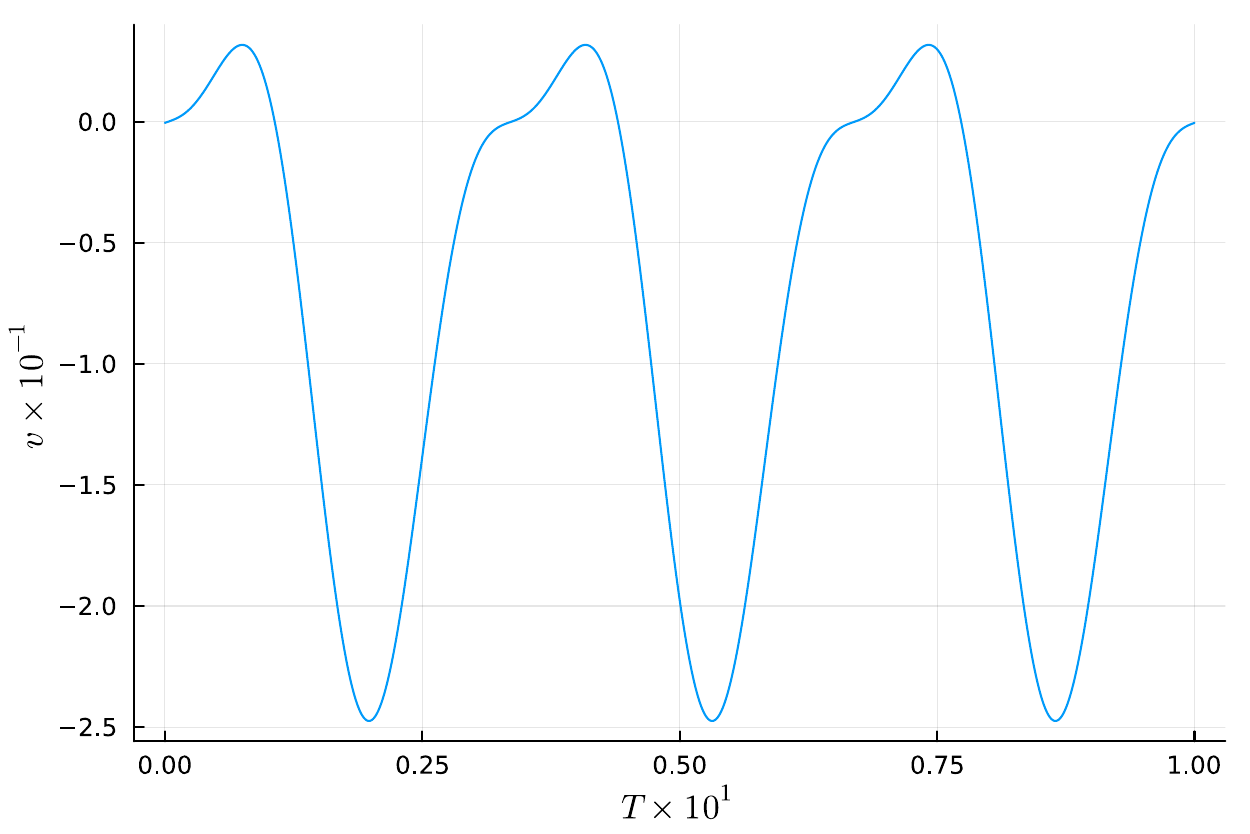}
    \includegraphics[width=0.3\textwidth]{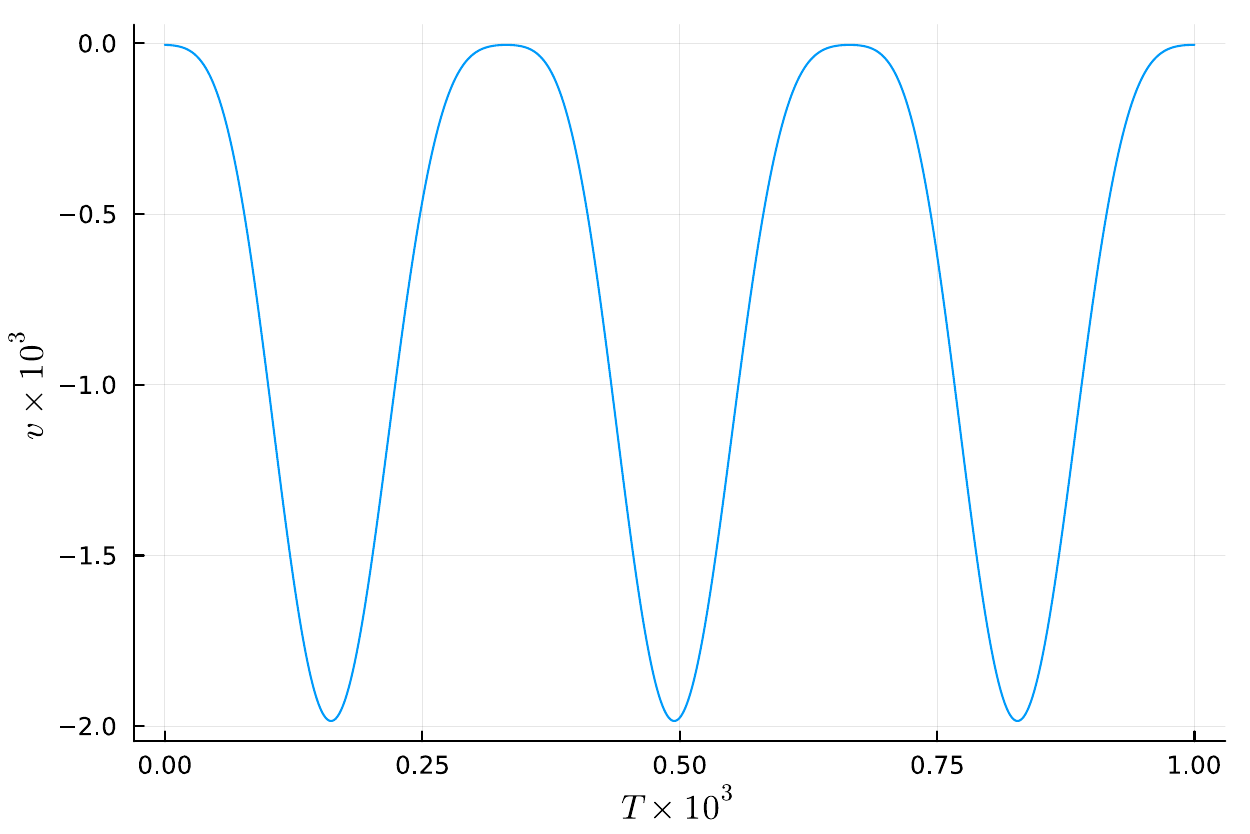}
    \caption{(Left) \(T = 10^{-1}\), (Center) \(T = 10\), (Right) \(T = 10^{3}\).}
    \label{fig-waveform-comparison}
\end{figure}

To explicitly handle periodicity, we introduce a convenient integer-valued function $\nu:\mathbb{R}\to\mathbb{Z}$ and an angular variable $\theta:\mathbb{R}\to[0,2\pi/\omega)$ by setting
\begin{align}
    \nu(x) &= \left\lfloor \frac{\omega x}{2\pi} \right\rfloor,\quad
    \theta(x) = x - \frac{2\pi}{\omega}\nu(x).
\end{align}
In what follows, we simply write $\nu = \nu(x)$, $\theta = \theta(x)$. Using these variables, the integral admits the following asymptotic expansions as $\varepsilon\to0$:
\begin{align}
	\int_0^x \frac{\md s}{\kappa_1+\kappa_0 \cos \omega s}
	&=\frac{2}{\omega\sqrt{\kappa_1^2-\kappa_0^2}}\left(\frac{\pi \nu}{2}+ \tan^{-1}\left(\sqrt{\frac{\kappa_1-\kappa_0}{\kappa_1+\kappa_0}}\tan\frac{\omega \theta}{2}\right)\right)\\
	&=\frac{x-\theta+\pi/\omega}{\sqrt{\kappa_1^2-\kappa_0^2}}+o(1).
\end{align}
Additionally, when $\theta=0$, we have
\begin{align}
	\int_0^x \frac{1}{\kappa_1+\kappa_0 \cos \omega s}\md s
	=\frac{x}{\sqrt{\kappa_1^2-\kappa_0^2}}.
\end{align}

Thus, we define a concise function $\ell:\mathbb{R}\to\mathbb{R}$ by
\begin{align}
	\ell(x)=
	\begin{cases}
		x,&\text{if }\theta=0,\\
		x-\theta+\pi/\omega,&\text{otherwise},
	\end{cases}
\end{align}
yielding a simple asymptotic form:
\begin{align}
	\int_0^x \frac{1}{\kappa_1+\kappa_0 \cos \omega s}\md s
	=\frac{\ell(x)}{\sqrt{\kappa_1^2-\kappa_0^2}} + o(1).
\end{align}

In the limit $\varepsilon\to 0$, we obtain the following approximation of the velocity $v$:
\begin{align}\label{eq-v-amp-diff-small-refined}
	v
	=-\frac{K(1-\cos \omega x)}{2\pi\omega^2}\left(-\omega \sin \omega x+\left(1-\frac{2(x-\ell(x))}{T} \right)(1-\cos \omega x)\right)+o(1).
\end{align}

Thus, when $|\kappa_1|-|\kappa_0|$ is sufficiently small, the velocity field $v$ can be decomposed into the following three functions:
\begin{align}
	P_1(x)&=- \sin \omega x(1-\cos \omega x)=-\sin \omega x+\frac{1}{2}\sin 2\omega x,\\
	P_2(x)&=(1-\cos \omega x)^2=\frac{3}{2} -2\cos\omega x+\frac{1}{2}\cos 2\omega x,\\
	A(x)&=1-\frac{2}{T}\left(x-\ell(x)\right).
\end{align}
Employing these functions, we have the following expression for $v$:
\begin{align}
	v=-\frac{K}{2\pi \omega^2}\left(-\omega P_1(x)+A(x)P_2(x)\right)+o(1),\quad\text{as }\varepsilon\to0.
\end{align}
Here, since the functions $P_1$ and $P_2$ include not only the fundamental oscillations $\cos \omega x$ and $\sin \omega x$, but also the double-angle terms $\cos 2\omega x$ and $\sin 2\omega x$ derived from their products, the response of the point vortex velocity to the curvature can be identified as nonlinear when $|\kappa_1|-|\kappa_0|$ is small. (see Figure~\ref{fig-waveform-comparison}).
Especially, the affine similarity law, which was applicable in the situation with sufficiently small $\kappa_0$ considered earlier, does not hold in general for the present scenario.
A noteworthy point is that the terms controlling $v$ vary according to the scale of $T$.
Specifically, when the period $T$ is small, $\omega P_1(x)$ dominates, whereas $A(x)P_2(x)$ becomes dominant as $T$ increases.
Furthermore, for each $x\in\re$, the function $A$ consistently fulfills $A(x)\in (1-1/N,1+1/N)$, independently of $T$.
Therefore, in the regime of large $T$ and large $N$, we have $A(x)\approx 1$ and the dominance of $P_2$, while for small values of $N$, the amplitude adjustment by $A$ becomes significant.

\begin{figure}[htbp]
    \centering
    \includegraphics[width=0.3\textwidth]{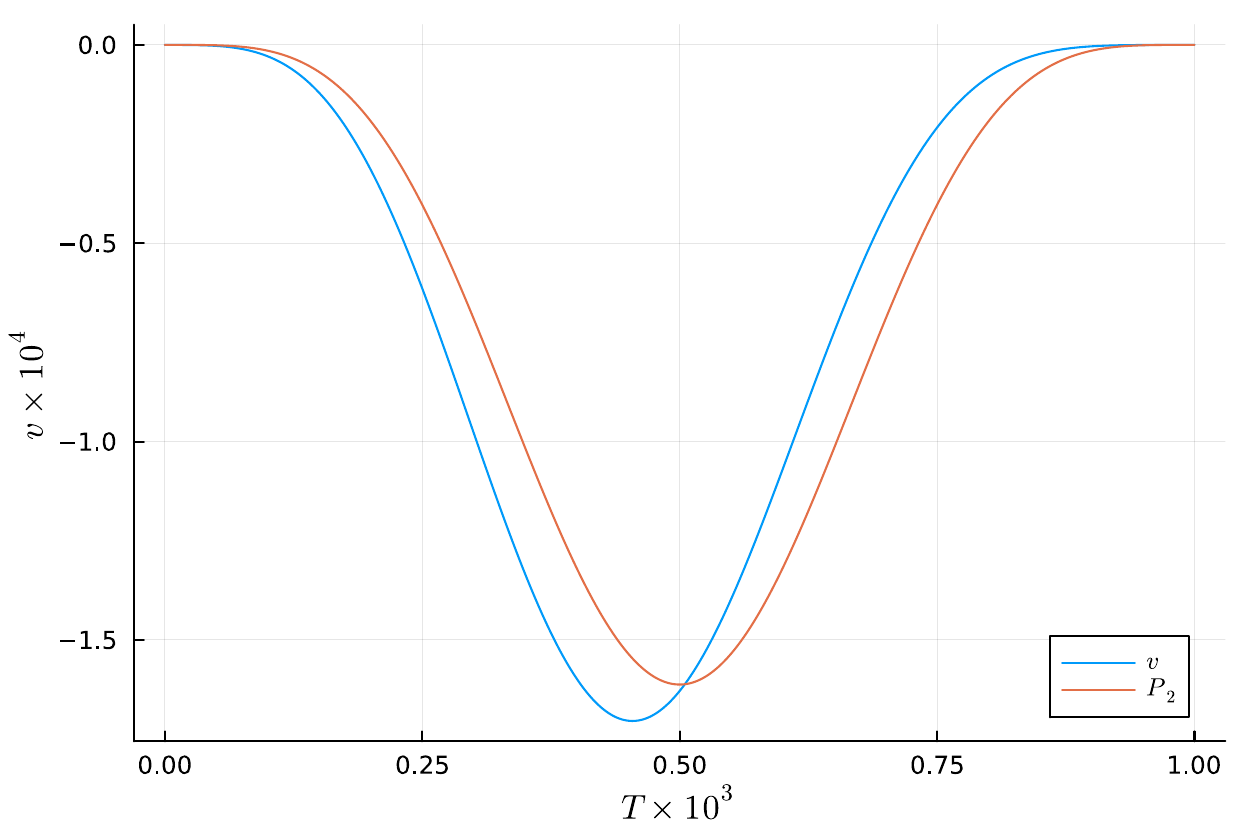}
    \caption{Comparison of the velocity field $v$ (blue) and the characteristic waveform $P_2$ (red) for $N=1$. The figure demonstrates the asymmetry induced by the modulation function $A$.}
    \label{fig-asymmetric-waveform-N1}
\end{figure}

In Figure~\ref{fig-asymmetric-waveform-N1}, we illustrate $v$ alongside $P_2$ for the case of $N=1$. 
This figure indicates that the minimum of $P_2$ is located at $x=\pi/\omega$, whereas the minimum point $x_*$ of $v$ is shifted to $x_*<\pi/\omega$. 
This occurs because the function $A$ nonuniformly modulates the amplitude of the $\pi/\omega$-periodic function $P_2$ depending on $x$, thereby inducing an asymmetry in the waveform.
In summary, the nonlinear effect of curvature on the point vortex velocity in the regime where $|\kappa_1|-|\kappa_0|$ is small is characterized by the two distinctive periodic functions $P_1$ and $A P_2$. Notably, the effect of waveform asymmetry induced by $A$ becomes significant especially when $N=1$ and $T$ is sufficiently large.


\section{Conclusion}
\label{sec-conclusion}
We have presented an analytic formula for the hydrodynamic Green function~$G$ and the associated Robin function~$R$ on every surface that admits a hydrodynamic Killing vector field.  Theorem~\ref{thm-Analytic-form} gives a unified representation of $G$ for fourteen canonical Riemann surfaces and, via conformal equivalence, for any surface in the same class.  These expressions automatically satisfy the slip boundary condition, the completeness of $-\jgrad G(\,\cdotp\,,x_{0})$, and the normalisation prescribed in Definition~\ref{def-HGF}.  In particular, the formula recovers the classical planar Biot--Savart kernel, extends it to compact and non-compact curved surfaces, and provides closed-form expressions for tori and annuli where explicit Green functions were previously unavailable.

As an immediate application, we have quantified the action of curvature on a single point vortex.  Equation~\eqref{eq-Robin-curvature} shows that, up to a surface-area constant, the vorticity of the Robin function coincides with the Gauss curvature.  Hence the vortex trajectory
\[
\dot q=-\jgrad R(q)
\]
is determined entirely by the geometric data of the surface.  Section~\ref{subsec-general-curvature} interprets this relation as a Hamiltonian flow driven by curvature, while Sections~\ref{subsec-linear-response} and~\ref{subsec-nonlinear-response} illustrate linear and nonlinear response regimes through explicit asymptotic expansions on a rectangular torus.  These computations clarify how the background curvature profile and lattice parameters influence vortex speed, waveform symmetry, and phase shift.

The analytic formula for $G$ supplies a concrete tool for further studies of Euler–Arnold flows on curved domains, including multi-vortex interactions, stability analysis, and numerical verification.  Because each term in the formula is expressed with elementary functions or standard special functions, the results can be implemented directly in symbolic or numerical software without additional simplification.

\vspace{0.5cm}
\noindent
{\bf Acknowledgments.}\ 
Research of YS was partially supported by Grant-in-Aid for JSPS Fellows 21J00025, Japan Society
for the Promotion of Science (JSPS).
\bibliographystyle{amsplain} 
\providecommand{\href}[2]{#2}

\end{document}